\documentclass[12pt,letter]{article}

\usepackage{amssymb,epsfig,latexsym}

\usepackage{graphicx,amsmath,amssymb,latexsym,epsfig}

\def\cH{{\cal H}}
\def\cB{{\cal B}}
\def\cX{{\cal X}}

\def\cK{{\cal K}}

\def\cC{{\cal C}}

\def\N{\mathbb{N}}
\def\R{\mathbb{R}}
\def\B{{\cal B}}
\def\cP{{\cal P}}

\def\and{\quad\mbox{and}\quad}
\def\ind{{\bf 1}}
\def\bn{{\bf n}}

\def\PP{{\mathrm P}}
\def\EE{{\mathrm E}}
\def\Z{{\mathbb Z}}
\def\var{{\mathrm{Var}}}

\newtheorem{corollary}{Corollary}
\newtheorem{theorem}{Theorem}
\newtheorem{lemma}{Lemma}

\def\ep{\hfill $\Box$}

\title{On the birth-and-assassination process, with an application to scotching a rumor in a network}
\author{Charles Bordenave\footnote{
Universit\'e de Toulouse \& CNRS - bordenave@math.univ-toulouse.fr}}
\begin{document}

\maketitle
\begin{abstract}
We give new formulas on the total number of born particles in the
stable birth-and-assassination process, and prove that it has a
heavy-tailed distribution. We also establish that this process is
a scaling limit of a process of rumor scotching in a network, and
is related to a predator-prey dynamics.


\noindent {\em MSC-class: } 60J80.

\end{abstract}

\section{Introduction}

\subsubsection*{Birth-and-assassination process}

The birth-and-assassination process was introduced by Aldous and
Krebs \cite{aldouskrebs}, it is a variant of the branching
process. The original motivation of the authors was then to
analyze a scaling limit of a queueing process with blocking which
appeared in database processing, see Tsitsiklis, Papadimitriou and
Humblet \cite{tsitsiklis}. In this paper, we show that the
birth-and-assassination process exhibits some heavy-tailed
distribution. For general references on heavy-tail distribution in
queueing processes, see for example Mitzenmacher
\cite{mitzenmacher} or Resnick \cite{resnick}. In this paper, we
will not discuss this application. Instead, we will show that the
birth-and-assassination process is also the scaling limit of a
rumor spreading model which is motivated by network epidemics and
dynamic data dissemination (see for example,  \cite{surveySIR},
\cite{andersson}, \cite{dynamicinformation}).

\begin{figure}[htb]
\begin{center}
\includegraphics[angle=0,height = 6cm]{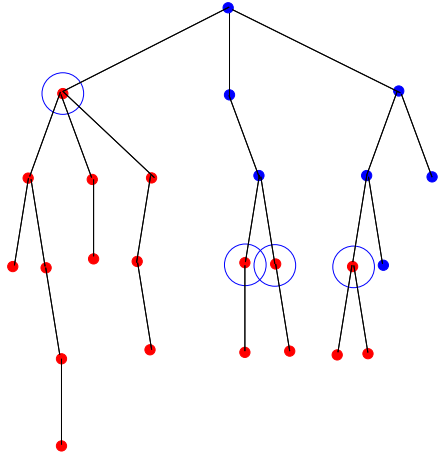}
\caption{Illustration of the birth-and-assassination process, living particles are in red, dead particles in blue, particles at risk are encircled.}
\end{center}
\end{figure}

We now reproduce the formal definition of the
birth-and-assassination process from \cite{aldouskrebs}.  Let
$\N^f = \cup_{k=0} ^{\infty} \N^k$ be the set of finite k-tuples
of positive integers (with $N^0 =\emptyset$). Let $\{\Phi_\bn\},
\bn \in \N^f$, be a family of independent Poisson processes with
common arrival rate $\lambda$. Let $\{K_\bn\}, \bn \in \N^f$, be a
family of independent, identically distributed (iid), strictly
positive random variables. Suppose the families $\{ \Phi_\bn \}$
and $\{ K_\bn \}$ are independent. The particle system starts at
time $0$ with only the ancestor particle, indexed by $\emptyset$.
This particle produces offspring at the arrival times of
$\Phi_\emptyset$, which enter the system with indices $(1)$,
$(2)$, $\cdots$ according to their birth order. Each new particle
$\bn$ entering the system immediately begins producing offspring
at the arrival times of $\{\Phi_\bn\}$, the offspring of $\bn$ are
indexed $(\bn, 1)$, $(\bn, 2)$, $\cdots$ also according to birth
order. The ancestor particle is {\em at risk} at time $0$. It
continues to produce offspring until time $D_\emptyset =
K_\emptyset$, when it dies. Let $k>0$ and let
$\bn=(n_1,\cdots,n_{k-1},n_k)$, $\bn'=(n_1,...,n_{k-1})$. When a
particle $\bn'$ dies (at time $D_{\bn'}$), $\bn$ then becomes at
risk; it continues to produce offspring until time $D_\bn =
D_{\bn'} + K_\bn$, when it dies. We will say that the
birth-and-assassination process is stable if with probability $1$
there exists some time $t < \infty$ with no living particle. The
process is unstable if it is not stable. Aldous and Krebs
\cite{aldouskrebs} proved the following:

\begin{theorem} [Aldous and Krebs]
\label{th:stab}Consider a birth-and-assassination process with offspring rate $\lambda$ whose killing distribution has
moment generating function $\phi$. Suppose $\phi$ is finite in some neighborhood of $0$. If $\min_{u >0}  \lambda u^{-1} \phi( u) < 1$ then
the process is stable. If $\min_{u >0}  \lambda u^{-1} \phi( u) > 1$ then the process is unstable.
\end{theorem}

The birth-and-assassination process is a variant the classical
branching process. Indeed, if instead the particle $\bn$ is at
risk not when its parent dies but when the particle $\bn$ was
born, then we obtain a well-studied type of branching process, refer to Athreya and Ney \cite{athreyaney}. The
populations in successive generations behave as the simple
Galton-Walton branching process with mean offspring equal to
$\lambda \EE K_\emptyset$, and so the process is stable if this
mean is less than $1$.  The birth-and-assassination process is a
variation in which the 'clock' which counts down the time until a
particle's death does not start ticking until the particle's
parent dies.

In this paper, we will pay attention to the special case where the killing distribution is an exponential distribution with intensity $\mu$. By a straightforward scaling argument,   a birth-and-assassination process with intensities $(\lambda,\mu)$ and a birth-and-assassination process with intensities $(\lambda \mu^{-1},1)$ where the time is accelerated by a factor $\mu$ have the same distribution. Therefore, without loss of generality, from now on, we will consider $\B$, a birth-and-assassination process with intensities $(\lambda, 1)$. As a corollary of Theorem \ref{th:stab}, we get
\begin{corollary}  [Aldous and Krebs]
\label{cor:stab} If $0 < \lambda < 1/4$,
the process $\B$ is stable. If $\lambda > 1/4$, the process $\B$ is unstable.
\end{corollary}

In the first part of this paper, we study the behavior of the
process $\B$ in the stable regime, especially as $\lambda$ get
close to $1/4$. We introduce a family of probability measures
$\{\PP_\lambda\}, \lambda >0,$ on our underlying probability space
such that under $ \PP_{\lambda}$, $\B$ is a
birth-and-assassination process with intensities $(\lambda, 1)$.
Let $\lambda \in (0,1/4)$, we define $N$ as the total number of
born particles in $\cB$ (including the ancestor particle) and
$$
\gamma (\lambda) = \sup \left\{ u \geq 0 :  \EE_{\lambda} N ^u < \infty \right \}.
$$
In particular, if $0< \gamma(\lambda) < \infty$, from Markov Inequality, for all $0 < \epsilon < \gamma(\lambda)$, there exists a constant $C \geq 1$ such that for all $t \geq 1$,
$$ \PP_{\lambda} ( N > t ) \leq C t^{-\gamma(\lambda) + \epsilon}.$$
The number $\gamma$ may thus be interpreted as a power tail
exponent. There is a simple expression for $ \gamma$.
\begin{theorem}
\label{th:powertail} For all $\lambda \in (0,1/4)$,
$$
\gamma (\lambda) = \frac{1 +  \sqrt {1 - 4 \lambda} } { 1-  \sqrt {1 - 4 \lambda} }.
$$
\end{theorem}
This result contrasts with the behavior of the classical branching
process, where for all $\lambda < 1$: there exists a constant
$c>0$ such that $\EE_\lambda \exp(c N ) < \infty$. This heavy tail
behavior of the birth-and-assassination process is thus a striking
feature of this process. Near criticality, as $\lambda  \uparrow
1/4$, we get $\gamma(\lambda) \sim 1$, whereas as $\lambda
\downarrow 0$, we find $\gamma(\lambda) \sim (2\lambda)^{-1}$. By
recursion, we will also compute the moments of $N$.
\begin{theorem}
\label{th:N}
\begin{itemize}
\item[(i)] For all $ p \geq 2$,  $\EE_{\lambda} N^p < \infty$
if and only if $ \lambda  \in (0,  p (p+1)^{-2})$.
 \item[(ii)] If
$ \lambda \in (0,1/4]$,
\begin{equation} \label{eq:N1}
\EE_{\lambda} N = \frac{2}{1 + \sqrt{1 - 4 \lambda}}.
\end{equation}
\item[(iii)]
If $\lambda \in (0,2/9)$,
\begin{equation} \label{eq:N2}
\EE_{\lambda} N^2 =\frac{2}{3 \sqrt{ 1 - 4 \lambda} -1}.
\end{equation}
\end{itemize}
\end{theorem}
Theorem \ref{th:N}(i) is consistent with Theorem
\ref{th:powertail}: $\lambda  \in (0,  p (p+1)^{-2})$ is
equivalent to $p \in [1 , (1 +  \sqrt {1 - 4 \lambda})(  1-  \sqrt
{1 - 4 \lambda})^{-1})$ . Theorem \ref{th:N}(ii) implies a
surprising discontinuity of the function $\lambda \mapsto
\EE_\lambda N $ at the critical intensity $ \lambda = 1/4$:
$\lim_{\lambda \uparrow 1/4} \EE_{\lambda} N = 2$. Again, this discontinuity contrasts with what happens
in a standard Galton-Watson process near criticality, where for
$0< \lambda < 1$, $\EE_\lambda N  = (1 - \lambda) ^{-1}$. We will
prove also that  this discontinuity is specific to $\lambda =1/4$
and for all $p \geq 2$, $\lim_{\lambda \uparrow p(p+1)^{-2}}
\EE_{\lambda} [ N^p] = \infty$. We will explain a method to
compute all integers moments of $N$ by recursion. The third moment
has already a complicated expression (see \S \ref{subsub:high}).
From Theorem \ref{th:N}(ii), we may fill the gap in Corollary
\ref{cor:stab}.
\begin{corollary}
\label{cor:N}
 If $\lambda = 1/4$,
the process $\B$ is stable.
\end{corollary}

In Section \ref{sec:birth-and-assassination}, we will prove
Theorems \ref{th:powertail} and \ref{th:N} by exhibiting a
Recursive Distributional Equation (RDE) for a random variable
related to $N$. Unfortunately, our method does not give much
insights on the heavy-tail phenomena involved in the
birth-and-assassination process.

\subsubsection*{Rumor scotching process}

We now define the rumor scotching process on a graph. It is a
nonstandard SIR dynamics (see for example \cite{surveySIR} or
\cite{andersson} for some background). This process represents the
dynamics of a rumor/epidemic spreading on the vertices of a graph
along its edges. A vertex  may be unaware of the rumor/susceptible
(S), aware of the rumor and spreading it as true/infected (I), or
aware of the rumor and trying to scotch it/recovered (R).

More formally, we fix a connected graph $G = (V,E)$, and let
$\cP_V$ denote the set of subsets of $V$ and $\cX=   ( \cP_V
\times \{ S, I , R\} )^V$. The spread of the rumor is described by
a Markov process on $\cX$. For $X = (X_v)_{v \in V}  \in \cX$,
with $X_v = (A_v, s_v)$, $A_v$ is interpreted as the set of
neighbors of $v$ which can change the opinion of $v$ on the
veracity of the rumor.  If $(uv) \in E$, we define the operations
$E_{uv}$ and $E_v$ on $\cX$  by $(X + E_{uv} )_w  = (X - E_{v}
)_w= X_w$, if $w \ne v$ and $(X + E_{uv} )_v =  (A_{v} \cup \{u\},
I)$,  $(X - E_v)_v =  (\emptyset, R)$. Let $\lambda >0$ be a fixed
intensity, the rumor scotching process is the Markov process with
generator:
\begin{eqnarray*}
K(X, X + E_{uv}) & =  & \lambda \ind ( s_u =  I ) \ind ( (u,v) \in E )  \ind ( s_v  \ne  R ), \\
 K(X, X - E_{v}) & =  & \ind ( s_v =  I ) \sum_{ u \in A_v } \ind ( s_u = R),
\end{eqnarray*}
and all other transitions have rate $0$. Typically, at time $0$, there is non-empty finite set of $I$-vertices and there is a vertex $v$ such that $A_v$ contains a $R$-vertex. The absorbing states of this process are the states without $I$-vertices.  The case when at time $0$, $A_v$ is the set to all neighbors of $v$ is interesting in its own (there, $A_v$ does not evolve before $s_v =R$).

If $G$ is the infinite $k$-ary tree this process has been analyzed
by Kordzakhia \cite{kordzakhia} and it was defined there as the
{\em chase-escape model}. It is thought as a predator-prey
dynamics: each vertex may be unoccupied (S), occupied by a prey
(I) or occupied by a predator (R). The preys spread on unoccupied
vertices and predators spread on vertices occupied by preys.  If
$G$ is the $\Z^d$-lattice and if there is no $R$-vertices, the
process is the original Richardson's model  \cite{richardson}.
With $R$-vertices, this process is a variant of the two-species
Richardson model with prey and predators, see for example
H{\"a}ggstr{\"o}m and Pemantle \cite{haggstrom98}, Kordzakhia  and
Lalley \cite{kordzakhia05}. Nothing is apparently known on this
process.

In Section \ref{sec:rumor}, we show that the
birth-and-assassination process is  the scaling limit, as $n$ goes
to infinity, of the rumor scotching process when $G$ is the
complete graph over $n$ vertices and the intensity is $\lambda /
n$ (Theorem \ref{th:rq}).

\section{Integral equations for the birth-and-assassination process} \label{sec:birth-and-assassination}

\subsection{Proof of Theorem \ref{th:N} for the first moment} \label{subsec:N1}
In this paragraph, we prove Theorem  \ref{th:N}(ii). Let $X(t) \in [0,+\infty]$ be the total
number of born particles in the process $\B$ given that the root
cannot die before time $t$, and $Y(t)$ be the total number of born
particles given that the root dies at time $t$. By definition, if
$D$ is an exponential variable with mean $1$ independent of $Y$,
then $N \stackrel{d}{=}  X(0) \stackrel{d}{=}  Y(D)$, where the
symbol $\stackrel{d}{=}$ stands for distributional equality. We
notice also that the memoryless property of the exponential
variable implies $X(t) \stackrel{d}{=} Y(t+D)$. The recursive
structure of the birth-and-assassination process leads to the
following equality in distribution
\begin{equation*}
\label{eq:RDE1}
Y(t) \stackrel{d}{=}  1 +  \sum_{i: \xi_i \leq t } X_i ( t - \xi_i) \stackrel{d}{=} 1 +  \sum_{i: \xi_i \leq t } X_i ( \xi_i),
\end{equation*}
where $\Phi = \{\xi_i\}_{i \in \N}$ is a Poisson point process of
intensity $\lambda$ and $(X_i), i \in \N,$ are  independent copies
of $X$. Note that since all variables are non-negative, there is
no issue with the case $Y(t) = + \infty$. We obtain the following
RDE for the random function $Y$:
\begin{equation}
\label{eq:RDE2}
Y(t) \stackrel{d}{=} 1 +  \sum_{i: \xi_i \leq t }  Y_i ( \xi_i + D_i),
\end{equation}
where  $Y_i,$ and $D_i$ are independent copies of $Y$ and $D$ respectively. This last RDE is the cornerstone of this work.

Assuming that $\EE_\lambda N < \infty$ we first prove that
necessarily $\lambda \in (0,1/4)$.  For convenience, we often drop
the parameter $\lambda$ in $\EE_{\lambda}$ and other objects
depending on $\lambda$. From Fubini's theorem, $\EE X(0) = \EE N =
\int_{0}^{\infty} \EE Y(t) e^{-t} dt$ and therefore $\EE Y(t) <
\infty$ for almost all $t \geq 0$. Note however that since $t
\mapsto Y(t)$ is monotone for the stochastic domination, it
implies that $\EE Y(t) < \infty$ for all $t > 0$. The same
argument gives the next lemma.
\begin{lemma}
\label{le:YNfinite} Let $t > 0$ and $u > 0$, if $\EE [ N^u]<
\infty$ then $\EE  [ Y(t)^u ]< \infty$.
\end{lemma}
Now, taking expectation in (\ref{eq:RDE2}), we get
\begin{eqnarray*}
\EE Y(t) & = & 1 + \lambda \int_0 ^ t \int_{0} ^ {\infty} \EE Y (x+s) e^{-s} ds dx.
\end{eqnarray*}
Let $f_1(t) = \EE Y(t)$, it satisfies the integral equation, for all $t \geq 0$,
\begin{equation}
\label{eq:g}
f_1(t) = 1+ \lambda \int_0 ^ t e^x \int_{x} ^ {\infty} f_1(s) e^{-s} ds dx.
\end{equation}
Taking the derivative once and multiplying by $e^{-t}$, we get:
$
f'_1 (t) e^{-t} = \lambda \int_t ^ {\infty} f_1(s) e^{-s} ds.
$
 Then, taking the derivative a second time and  multiplying  by $e^{t}$:
$
f_1'' (t) - f_1' (t) =   - \lambda f_1(t).
$
So, finally, $f_1$ solves a linear ordinary differential equation of the second order
\begin{equation}
\label{eq:ODE1}
x'' - x' + \lambda x = 0,
\end{equation}
with initial condition $x(0) = 1$. If $\lambda > 1/4$ the
solutions of (\ref{eq:ODE1}) are  $$x(t) =  e^{t/2} (  \cos ( t
\sqrt{4\lambda -1 } )  + a \sin (t \sqrt{4\lambda -1 } )), $$
for some constant $a$. Since $f_1 (t)$ is necessarily positive,
this leads to a contradiction and $\EE N = \infty$. Assume now that $0 < \lambda < 1/ 4$ and let
\begin{equation}
\label{eq:alphabeta}
\Delta =  \sqrt { 1 - 4 \lambda }  \; \hbox{ , } \quad  \alpha = \frac{ 1 - \Delta}{ 2}   \quad  \hbox{ and } \quad \beta = \frac{ 1 +\Delta}{  2}.
\end{equation}
$(\alpha,\beta)$ are the roots of the polynomial $X^2 - X + \lambda = 0$.  The solutions of (\ref{eq:ODE1}) are
$$
x_a ( t )  = (1-a) e^{\alpha t } + a e^{\beta t}
$$
for some constant $a$. Whereas, for  $\lambda = 1/4$, $\alpha = 1/2$ and the solutions of  (\ref{eq:ODE1}) are
$$
x_a ( t )  = (at +1 ) e^{t /2 }.
$$
For $0 < \lambda \leq 1/4$, we check easily that the functions $x_a$ with $a \geq 0$ are the nonnegative solutions of the integral equation (\ref{eq:g}).

It remains to prove that if $0 < \lambda \leq 1/ 4$ then $\EE N < \infty$ and
$
f_{1}  ( t ) =  e^{\alpha t } .
$
Indeed, then $\EE N = \int_0^ {\infty} f_1(t) e^{-t} dt =  (1 - \alpha)^{-1} $ as stated in Theorem \ref{th:N}(ii). To this end, define $f_1^{(n)} (t) = \EE \min( Y(t), n)$, from (\ref{eq:RDE2}),
$$
\min ( Y(t),n) \leq_{st} 1 +  \sum_{i: \xi_i \leq t }  \min(Y_i ( \xi_i + D_i),n).
$$
Taking expectation, we obtain, for all $t \geq0$,
\begin{equation}
\label{eq:gn}
f_1^{(n)}(t) \leq 1 + \lambda \int_0 ^ t e^x \int_{x} ^ {\infty} f_1^{(n)}(s) e^{-s} ds dx.
\end{equation}
We now state a lemma which will be used multiple times in this paper. We define
\begin{equation}
\label{eq:gamma} \overline \gamma (\lambda) = (1 + \Delta )  /  (1
- \Delta )= \beta / \alpha.
\end{equation}
Let  $1 < u < \overline \gamma$ (or equivalently $\lambda < u
(u+1)^{-2} $),  we define $\cH_u$, the set of measurable functions
$h : [0,\infty) \to [0,\infty)$ such that $h$ is non-decreasing
and $\sup_{t \geq 0} h(t) e^{-u\alpha t} < \infty$. Let $C> 0$, we
define the mapping from $\cH_u$ to $\cH_u$,
$$
\Psi : h \mapsto   C e^{u \alpha t }    + \lambda \int_0 ^ t e^x \int_x ^ {\infty}h (s) e^{-s} ds dx.
$$
In order to check that $\Psi$ is indeed a mapping from $\cH_u$ to
$\cH_u$, we use the fact that if $1 < u <\overline \gamma$, then
$u \alpha < 1$. Note also that if $1 < u <\overline \gamma$, then
$u\alpha - \lambda - u^2 \alpha^2>0$. If $\lambda =1/4$, we also
define the mapping from $\cH_1$ to $\cH_1$,
$$
\Phi : h \mapsto  1    +  \frac 1 4 \int_0 ^ t e^x \int_x ^ {\infty}h (s) e^{-s} ds dx.
$$
(recall that for $\lambda =1/4$, $\alpha = 1/2$).

\begin{lemma}
\label{le:powerupp2}\begin{enumerate}
\item[(i)]
Let  $1 < u < \overline \gamma$ and $f \in \cH_u$ such that $f \leq \Psi(f)$. Then for all $t \geq 0$,
$$
f(t) \leq  C \frac{ u \alpha ( 1- u \alpha) } { u\alpha - \lambda - u^2 \alpha^2} e ^{u \alpha t}   -C \frac{\lambda} { u\alpha - \lambda - u^2 \alpha^2}  e^{\alpha t} .
$$
\item[(ii)] If $\lambda =1/4$ and $f \in \cH_1$ is such that $f
\leq \Phi(f)$, then for all $t \geq 0$,
$$
f(t) \leq  e ^{t/2} .
$$
\end{enumerate}
\end{lemma}
Before proving Lemma \ref{le:powerupp2}, we conclude the proof of
Theorem \ref{th:N}$(ii)$. For $0 < \lambda < 1/4$, from
(\ref{eq:gn}), we may apply Lemma \ref{le:powerupp2}(i)  applied
to $1 < u < \beta /\alpha$, $C=1$. We get that $$f_1^{(n)} (t)
\leq  C_u e^{ \alpha u t}$$ for some $C_u >0$. The monotone
convergence theorem implies that $f_1 (t) = \lim_{n \to \infty}
f_1^{(n)} (t)$ exists and is bounded by $C_u e^{ \alpha u t}$.
Therefore $f_1$ solves the integral equation (\ref{eq:g}) and is
equal to $x_a$ for some $a \geq 0$. From what precedes, we get
$x_a(t) \leq C_u e^{ \alpha u t}$, however, since $ \alpha u <
\beta$, the only possibility is $a =0$ and $f_1 (t) = e^{\alpha t
}$.

Similarly, if  $\lambda = 1/4$, from Lemma \ref{le:powerupp2}$(ii)$,
$
f_1 (t) \leq e^{t/2}.
$
This proves that $f_1$ is finite, and we thus have $f_1 = x_a $ for some $a \geq 0$. Again, the only possibility is $a=0$ since $x_a (t) \leq e^{t/2}$ implies $a = 0$.

\noindent{{\em Proof of Lemma \ref{le:powerupp2}.} $(i)$.  The
fixed points of the mapping $\Psi$ are the functions $h_{a,b}$
such that
$$
h_{a,b} (t) = a e^{\alpha t} + b e^{\beta t } + C \frac{u \alpha ( 1  - u \alpha) }{ u\alpha - \lambda - u^2 \alpha^2} e ^{u \alpha t},
$$
with $a + b + C \frac{u \alpha ( 1  - u \alpha) }{   u\alpha -
\lambda - u^2 \alpha^2} = C$. The only fixed point in $\cH_u$ is
$h_*  := h_{a_*,0}$ with $a_* =  - C  \lambda /(   u\alpha -
\lambda - u^2 \alpha^2)$. Let $ \cC_u$ denote the set of
continuous functions in $\cH_u$, note that $\Psi$ is also a mapping
from $ \cC_u$ to $\cC_u$. Now let  $g_0 \in \cC_u$ and for $k \geq
1$, $g_{k} = \Psi (g_{k-1})$. We first prove that for all $t\geq0$
, $\lim_k g_k (t) = h_*(t)$. If $1 < u < \overline \gamma$ then $u
\alpha ( 1- u \alpha) > \lambda$ and $\frac{u \alpha ( 1  - u
\alpha) }{   u\alpha - \lambda - u^2 \alpha^2}$ is positive. We
deduce easily that if $ g_0 (t) \leq L e^{u \alpha t} $ then $g_1
(t) = \Psi (g) (t)  \leq  C e^{u \alpha t } +  \frac{ L \lambda}{
u \alpha ( 1- u \alpha) }  (e^{u \alpha t } -1) \leq L_1  e^{u
\alpha t} $, with $L_1 =  (C +  \frac{ L \lambda}{ u \alpha ( 1- u
\alpha) } )$. By recursion, we obtain that $\limsup_k g_k (t) \leq
L_{\infty} e^{u\alpha t}$, with $L_{\infty} = C u\alpha (1 - u
\alpha) / (u\alpha  - \lambda - u^2 \alpha^2)  < \infty$. From
Arzela-Ascoli's theorem, $(g_k)_{k \in \N}$ is relatively compact in
$\cC_u$ and any accumulation point  converges to $h_*$ (since
$h_*$ is the only fixed point of $\Psi$ in $ \cC_u$).

Now since $f \in \cH_u$, there exists a constant $ L  >0$ such
that for all $t\geq 0$,  $ f (t) \leq g_0 (t) := L e^{u \alpha
t}$. The monotonicity of the mapping $\Psi$ implies that $\Psi(f)
\leq \Psi(g_0) = g_1$. By assumption, $f \leq \Psi(f)$ thus by
recursion $f \leq \lim_n g_n = h_*$.

\noindent $(ii)$. The function $x_0  ( t ) =  e^{t /2} $ is the
only fixed point of $\Phi$ in $\cH_1$. Moreover, if $ g (t) \leq C
e^{t/2}$ then we also have $\Phi(g) (t) \leq  C e^{t/2}$.  Then,
if $g$ is continuous, arguing as above, from Arzela-Ascoli's
theorem, $(\Phi^k(g))_{k \in \N}$ converges to $x_0$. We conclude
as in (i). \ep

\subsection{Proof of Theorem \ref{th:N}(i)}

We define $f_p (t) = \EE_\lambda [ Y(t)^p ]$. As above, we often
drop the parameter $\lambda$ in $\EE_{\lambda}$ and  other objects
depending on $\lambda$.
\begin{lemma}
\label{le:rec}
Let $p \geq 2$, there exists a polynomial $Q_p$ with degree $p$ such that for all $t >0$,
\begin{itemize}
\item[(i)]
If $\lambda \in (0,p(p+1)^{-2})$, then $f_p (t) = Q_p ( e^{\alpha t})$.
\item[(ii)]
If $\lambda \geq p(p+1)^{-2}$, then $f_p (t) = \infty$,
\end{itemize}
\end{lemma}
Note that if such polynomial $Q_p$ exists  then $Q_p(x) \geq 1$ for all $x \geq 1$. Note also that $\lambda \in (0,p(p+1)^{-2})$ implies that $p < \overline \gamma = \beta / \alpha$ (where $\overline \gamma$ was defined by (\ref{eq:gamma})),  and thus $p \alpha < \beta < 1$.  Hence Lemma \ref{le:rec} implies Theorem \ref{th:N}(i) since $\EE [N^p]  = \int f_p (t) e^{-t} dt$.

Let $\kappa_p (X)$ denote the $p^{th}$ cumulant of a random
variable $X$ whose moment generating function is defined in a
neighborhood of $0$: $\ln \EE e^{ \theta X} = \sum_{p \geq 0}
\kappa_p (X) \theta^p / p!$. In particular $\kappa_0 (X) = 0$,
$\kappa_1 (X) = \EE X$ and $\kappa_2 (X) = \var X$.  Using the
exponential formula
\begin{equation}\label{eq:expformula} \EE \exp {\sum_{\xi_i \in \Phi} h(\xi_i,
Z_i) } = \exp ( \lambda \int_0 ^{\infty} (\EE e^{h(x,Z)} -1) dx),
\end{equation}
 valid for all non-negative function $h$ and iid variables
$(Z_i), i \in \N$, independent of $\Phi = \{\xi_i\}_{i \in \N}$ a
Poisson point process of intensity $\lambda$, we obtain that for
all $p \geq 1$,
\begin{equation}\label{eq:cumulPoi}
\kappa_p \left(\sum_{i : \xi_i \leq t} h (\xi_i,Z_i) \right) = \lambda \int_0 ^ t \EE h ^p (x,Z) dx.
\end{equation}
Due to this last formula, it will be easier to deal with the cumulant $g_p(t) = \kappa_p ( Y(t))$. By recursion, we will prove the next lemma which implies Lemma \ref{le:rec}.
\begin{lemma}
\label{le:rec2}
Let $p \geq 2$, there exists a polynomial $R_p$ with degree $p$, positive on $[1,\infty)$  such that, for all $t > 0$,
\begin{itemize}
\item[(i)] If $\lambda \in (0,p(p+1)^{-2})$, then $f_p (t)  <
\infty$ and $g_p (t) = R_p (e^{\alpha t})$. \item[(ii)] If
$\lambda \geq p(p+1)^{-2}$, then $f_p (t) = \infty$,
\end{itemize}
\end{lemma}

\noindent{{\em Proof of Lemma \ref{le:rec2}.}
In \S \ref{subsec:N1}, we have computed $f_p$ for $p=1$ and found $R_1 (x) = x$. Let $p \geq  2$ and assume now that the statement of the Lemma \ref{le:rec2} holds for $q =1, \cdots, p-1$. We assume first that $f_p (t) < \infty$, we shall prove that necessarily $\lambda \in (0,p(p+1)^{-2})$ and $g_p (t) = R_p (e^{\alpha t})$.  Without loss of generality we assume that $0 < \lambda < 1/4$. From Fubini's theorem, using the linearity of cumulants in (\ref{eq:RDE2}) and (\ref{eq:cumulPoi}), we get
\begin{eqnarray}
g_p (t) & = &  \lambda \int_0 ^ t \int_0 ^ {\infty} \EE [ Y ( x + s) ^p ] e^{-s} ds  dx \nonumber \\
& = &  \lambda \int_0 ^ t e^x \int_x ^ {\infty} f_p ( s)  e^{-s} ds  dx, \label{eq:gp}
\end{eqnarray}
(note that Fubini's Theorem implies the existence of $f_p(s)$ for all $s>0$). From Jensen inequality $ f_p (t) \geq g_1 (t) ^p = e^{ p \alpha t}$ and the integral $\int_x ^ {\infty}e^{ p \alpha s}  e^{-s} ds  dx$ is finite if and only if $p \alpha < 1 $. We may thus assume that $p  \alpha < 1$. We now recall the identity:
$ \EE X^p = \sum_{\pi} \prod_{I \in \pi} \kappa_{|I|}(X)$, where the sum is over all set partitions of $\{1,\cdots, p\}$, $I \in \pi$ means $I$ is one of the subsets into which the set is partitioned, and $|I|$ is the cardinal of $I$. This formula implies that $ \EE X^p = \kappa_p (X) + \Sigma_{p-1} (\kappa_1(X), \cdots, \kappa_{p-1}(X))$, where $\Sigma_{p-1}(x_1,\cdots,x_{p-1})$ is a polynomial in $p-1$ variables with non-negative coefficients and each of its monomial $\prod_{\ell = 1}^ k  x^{n_\ell}_{i_\ell}$ satisfies $\sum_\ell n_\ell i_\ell = p$.  Using the recurence hypothesis, we deduce from (\ref{eq:gp}) that there exists a polynomial $\tilde R_p (x)= \sum_{k=1} ^p r_k x^k$ of degree $p$ with $r_p >0$ such that
\begin{eqnarray}
g_p (t) & = &   \lambda \int_0 ^ t e^x \int_x ^ {\infty} \left(  g_p ( s)  e^{-s}  + \tilde R_p ( e^{\alpha s})  e^{-s} \right) ds  dx \nonumber \\
& = & \sum_{k=1} ^p \frac{ \lambda r_k}{k\alpha (1 - k \alpha)} e^{k \alpha t}  +   \lambda \int_0 ^ t e^x \int_x ^ {\infty} g_p ( s)  e^{-s} ds  dx \label{eq:gp2},
\end{eqnarray}
(recall that $ p \alpha < 1$). Now we take the derivative of this last expression, multiply by $e^{-t}$ and take the derivative again. We get that $g_p$ is a solution of the differential equation:
\begin{equation}\label{eq:ODEp}
x'' - x' + \lambda x =  - \sum_{k = 1}^p \lambda r_ {k} e^{k \alpha t},
\end{equation}
with initial condition $x(0) = 0$. Thus necessarily $g_p (t) = a
e^{\alpha t} + b e ^{\beta t } +  \varphi (t)$, where $\varphi(t)$
is a particular solution of the differential equation
(\ref{eq:ODEp}). Assume first that $\lambda \ne p (p+1)^{-2}$,
then it is easy to check that $ (p+1) \lambda - p \alpha$ and $ p
(p+1)^{-2} - \lambda$ are different from $0$ and have the same
sign. Looking for a function $\varphi$ of the form $\varphi(t) =
\sum_{k=1}^p c_k e^{k \alpha t}$ gives $c_k = - \lambda r_k ( k^2
\alpha^2 - k \alpha + \lambda)^{-1} =   \lambda r_k (k-1)^{-1} (
(k+1) \lambda - k \alpha)^{-1} $. If $\lambda >  p (p+1)^{-2}$
then $p \alpha > \beta$ and the leading term in $g_p$ is $c_p e^{p
\alpha t}$. However, if $\lambda >  p (p+1)^{-2}$, $c_p <0$ and
thus $g_p (t) < 0$ for $t$ large enough. This is a contradiction
with Equation (\ref{eq:gp}) which asserts that $g_p(t)$ is
positive.

We now check that if $0< \lambda < p (p+1)^{-2}$ then $f_p (t)$ is finite. We define $f_p^{(n)} (t) = \EE [ \min (Y(t), n)^p ]$. We use the following identity,
$$
\left( \sum_{i=1} ^N y_i  \right) ^ p =  \sum_{ i = 1} ^N \sum_{k=0} ^{p-1} { p-1 \choose k } y_i^{k+1}\left( \sum_{j \neq i} ^N y_i  \right) ^ {p-k-1}. $$
Then from (\ref{eq:RDE2}) we get,
\begin{align}
& (Y (t) - 1) ^p   \stackrel{d}{=} \label{eq:dev} \\
 & \quad  \sum_{ \xi_i \leq t }
Y_i ( \xi_i + D_i )^p  + \sum_{ \xi_i \leq t }   \sum_{k=0} ^{p-2}
{ p-1 \choose k } Y_i(\xi_i + D_i)^{k+1} \left( \sum_{\xi_j \neq
\xi_i \leq t } Y_j ( \xi_j + D_j )  \right) ^ {p-k-1}. \nonumber
\end{align} The recursion hypothesis implies that
there exists a constant $C$ such that $f_k (t) = Q_k (e^{\alpha t
}) \leq C e^{k\alpha t}$ for all $1 \leq k \leq p-1$. Thus, the
identity $Y(t)^p = (Y(t) - 1) ^p - \sum_{k=0}^{p-1} { p \choose k}
(-1) ^{p-k} Y(t)^k$ gives
\begin{eqnarray*}
f_p^{(n)} (t) & \leq & \EE [\min (Y(t) - 1,n)^p ] +  \sum_{k=0}^{p-1} { p \choose k} C e^{k\alpha t} \\
& \leq & \EE [\min (Y(t) - 1,n)^p ] +  C_1 e ^{p \alpha t}.
\end{eqnarray*}
From the recursion hypothesis, if $1 \leq k \leq p-1$, $$\int_0 ^
t \EE[ Y(x+D) ^ k] dx = \int_0 ^t e^{x} \int_x ^\infty f_k (s)
e^{-s} ds dx =  \tilde Q_k (e^{\alpha t}) \leq C e^{k \alpha t}$$
for some constant $C>0$. We take the expectation in (\ref{eq:dev})
and use Slyvniak's theorem to obtain
\begin{eqnarray*}
f_p^{(n)} (t)   & \leq  &   C_1 e ^{p \alpha t} +  \lambda \int_0 ^ {t} e^x \int_x^\infty f_p^{(n)}  (s) e^{-s} ds dx \\
&   & \;  + \; \lambda \int_0^t   \sum_{k=0} ^{p-2} { p-1 \choose k } \EE [ Y_i(x + D_i)^{k+1}] \EE  \left[ \Bigm( \sum_{\xi_j \leq t } Y_j ( \xi_j + D_j )   \Bigm)^ {p-k-1} \right]dx \\
&  \leq &   C_1 e ^{p \alpha t} +  \lambda \int_0 ^ {t} e^x \int_x^\infty f_p^{(n)}  (s) e^{-s} ds dx \\
&   & \;  + \; \lambda \sum_{k=0} ^{p-2} { p-1 \choose k } \tilde Q_{k+1} (e^{\alpha t}) \EE  [ (Y(t) -1  ) ^ {p-k-1} ] \\
&   \leq &    C_2 e ^{p \alpha t} +
\lambda \int_0 ^ {t} e^x \int_x^\infty f_p^{(n)}  (s) e^{-s} ds dx
\end{eqnarray*}
So finally for a suitable choice of $C$,
\begin{equation}\label{eq:boundfpn}|
f_p^{(n)} (t) \leq C e^{p  \alpha t}  + \lambda \int_0 ^ {t} e^x \int_x^\infty f_p^{(n)}  (s) e^{-s} ds dx.
\end{equation}
From Lemma \ref{le:powerupp2}, $f_p ^{(n)} (t) \leq C' e^{p \alpha t}$, and, by the monotone convergence theorem, $g_p (t) \leq f_p(t) \leq  C' e^{p \alpha t}$. From what precedes:
$g_p (t) = a e^{\alpha t} + b e ^{\beta t } +  \varphi (t)$, with $\varphi(t) = \sum_{k=1}^p c_k e^{k \alpha t}$, with $c_p >0$. If $b > 0$, since $\lambda >  p (p+1)^{-2}$ then $p \alpha < \beta$ and the leading term in $g_p$ is $b e^{\beta t}$ which is in contradiction with  $g_p (t)  \leq  C' e^{p \alpha t}$. If $b < 0$, this is a contraction with Equation (\ref{eq:gp}) which asserts that $g_p(t)$ is positive. Therefore $b =0$ and $g_p (t) = a e^{\alpha t} +  \varphi (t) = R_p (e^{\alpha t})$.

It remains to check that if $ \lambda = p (p+1)^{-2}$ then for all $ t >0$, $f_p (t) = \infty$. We have proved that, for all $\lambda < p (p+1)^{-2}$, $g_p (t) =   u_p(\lambda)  (p-1)^{-1} ( (p+1) \lambda - p \alpha)^{-1} e ^{p \alpha t} + S_{p-1} (e^{\alpha t})$, where $S_{p-1}$ is a polynomial of degree at most $p -1$ and $u_p (\lambda) >0$. Note that $\lim_{\lambda \uparrow p (p+1)^{-2}} (p+1) \lambda - p \alpha = 0$. A closer look at the recursion shows also that $u_p (\lambda)$ is a sum of products of terms in $\lambda$ and $\lambda (\ell-1) ^{-1}( (\ell+1) \lambda - \ell \alpha)^{ -1}$, with $2 \leq \ell \leq p-1$. In particular, we deduce that  $\lim_{\lambda \uparrow p (p+1)^{-2}} u_p (\lambda)  > 0$. Similarly, the coefficients of $S_{p-1}$ are equal to sums of products of integers and terms in $\lambda$ and $\lambda (\ell-1) ^{-1}( (\ell+1) \lambda - \ell \alpha)^{ -1}$, with $2 \leq \ell \leq p-1$. Thus they stay bounded as $\lambda$ goes to $p (p+1)^{-2}$ and we obtain, for all $t>0$,
\begin{equation}\label{eq:gpcritic}
\liminf_{\lambda \uparrow p (p+1)^{-2}}  f_p (t) \geq \lim_{\lambda \uparrow p (p+1)^{-2}}  g_p (t)  = \infty.
\end{equation}
Now, for all $t > 0$, the random variable $Y (t)$ is stochastically non-decreasing with $\lambda$. Therefore $\EE_\lambda [Y(t)^p]$ is non-decreasing and (\ref{eq:gpcritic}) implies that  $\EE_{1/4} [Y(t)^p] = \infty$. The proof of the recursion is complete.

\ep

\subsection{Proof of Theorem \ref{th:N}(iii)}
In this paragraph, we prove Theorem  \ref{th:N}(iii). Let $ \lambda \in (0,2/9)$, recall that $f_2 (t) = \EE Y(t)^2$ and $g_2 (t) = \var (Y(t))$. From (\ref{eq:gp}}) applied to $p=2$,
\begin{eqnarray*}
g_2 (t) & =  &  \lambda \int_0 ^ t \int_0 ^ {\infty} g_2 ( x + s)  e^{-s} + f_1^2 (x+s)  e^{-s} ds  dx .
\end{eqnarray*}
Since $f_1(t) = e^{\alpha t}$ and $\alpha^2 = \alpha - \lambda$,  $g_2$ satisfies the integral equation:
\begin{equation*}
\label{eq:f}
g_2( t) =  \frac{\lambda}{2(2\lambda - \alpha)}\left( e^{2 \alpha t }  -1\right) +   \lambda \int_0 ^ t e^{x} \int_x ^ {\infty} g_2 ( s)  e^{-s} ds  dx .
\end{equation*}
We deduce that $g_2$ solves an ordinary differential equation:
\begin{equation*}
\label{eq:ODE2}
x'' - x' + \lambda x = - \lambda e^{2\alpha t},
\end{equation*}
with initial condition $x(0) = 0$. Thus $g_2$ is of the form:
$
g_2(t) =  a e^{\alpha t} + b e^{\beta t} +\frac{\lambda}{3 \lambda -2 \alpha} e^{2 \alpha t}.
$
with $
a  + b  +\frac{\lambda}{3 \lambda -2 \alpha}  = 0
$. From Lemma \ref{le:rec2},  $b=0$ so finally
$$
g_2(t)  = \frac{\lambda}{3 \lambda - 2 \alpha}  \left( e^{2 \alpha t}  - e^{\alpha t} \right) \quad {and } \quad f_2(t)  = 2 \frac{2 \lambda -\alpha}{3 \lambda - 2 \alpha} e^{2 \alpha t} - \frac{\lambda}{3 \lambda - 2 \alpha}e^{\alpha t}.
$$
We conclude by computing $\EE N^2 = \int e^{-t} f_2 (t) dt$.

\subsection{Proof of Theorem \ref{th:powertail}}

As usual we drop the parameter $\lambda$ in $\EE_\lambda$. From (\ref{eq:gamma}), we have $\overline \gamma ( \lambda) = \frac{1 - 2 \lambda + \sqrt {1 - 4 \lambda} } {2 \lambda}$.  To prove Theorem \ref{th:powertail}, we shall prove two statements

\begin{equation}
\label{eq:power1}
\hbox{If $ \EE [ N^u]  < \infty$  then }  u \leq \overline \gamma,
\end{equation}

\begin{equation}
\label{eq:power2}
\hbox{If $1 \leq u < \overline \gamma$ then } \EE [ N^u]  < \infty.
\end{equation}

\subsubsection{Proof of (\ref{eq:power1}).}
Let $u \geq 1$, we assume that $ \EE[  N^u]Ê < \infty$. From Lemma \ref{le:YNfinite} and (\ref{eq:RDE2}), we get
 \begin{equation*}
\label{eq:RDEp}
\EE [ Y(t)^u  ] = \EE \left(1 +  \sum_{i: \xi_i \leq t }  Y_i ( \xi_i + D_i)\right)^u .
\end{equation*}
Let $f_u (t) = \EE [ Y(t)^u ]$. Taking expectation and using the inequality $(x + y)^u \geq x^u + y^u$, for all positive $x$ and $y$, we get:
\begin{eqnarray}
f_u (t)  &  \geq & 1 + \lambda \int_0 ^ t \EE f_u ( x + D ) dx \nonumber \\
& \geq & 1 + \lambda \int_0 ^ t e^x \int_x ^ {\infty} f_u (s) e^{-s} ds dx \label{eq:gplow}.
\end{eqnarray}
From Jensen's Inequality, $f_u (t) \geq f_1 (t) ^ u = e^{u \alpha t}$. Note that the integral $\int_x ^ {\infty} e^{\alpha u s}  e^{-s} ds$ is finite if and only if  $u <   \alpha^{-1}$.  Suppose now that $\overline \gamma < u < \alpha^{-1}$.  We use the fact: if  $ u > \overline \gamma$ then   $u^2 \alpha^2 - u \alpha + \lambda > 0$,
to deduce that there exists $0 < \epsilon < \lambda$ such that
\begin{equation}
\label{eq:defeps}
u^2 \alpha^2 - u \alpha + \lambda > \epsilon.
\end{equation}
Let $\tilde \lambda = \lambda - \epsilon$, $\tilde \alpha = \alpha (\tilde \lambda)$,  $\tilde \beta = \beta (\tilde \lambda)$,  we may assume that $\epsilon$ is small enough to ensure also that
\begin{equation}
\label{eq:deftilde2}
u \alpha > \tilde \beta.
\end{equation}
(Indeed, for all $\lambda \in (0,1/4)$, $\alpha (\lambda) \overline \gamma ( \lambda ) = \beta(\lambda)$ and the mapping  $\lambda \mapsto \beta (\lambda)$ is obviously continuous). We compute a lower bound from (\ref{eq:gplow}) as follows:
\begin{eqnarray}
f_u (t)  &  \geq & 1 + \tilde \lambda \int_0 ^ t e^x \int_x ^ {\infty} f_u (s) e^{-s} ds dx + \epsilon  \int_0 ^ t e^x \int_x ^ {\infty} f_u (s) e^{-s} ds dx \nonumber \\
& \geq & 1 + \tilde \lambda \int_0 ^ t e^x \int_x ^ {\infty} f_u (s) e^{-s} ds dx + \epsilon   \int_0 ^ t e^x \int_x ^ {\infty} e^{u \alpha s} e^{-s} ds dx   \nonumber \\
& \geq & 1 + C  ( e^{u \alpha t } -1 )  + \tilde \lambda \int_0 ^ t e^x \int_x ^ {\infty} f_u (s) e^{-s} ds dx \label{eq:gplow2},
\end{eqnarray}
with $C =   \epsilon  (u\alpha (1 - u\alpha))^{-1} >0$. We consider the mapping $\Psi : h \mapsto 1 + C  ( e^{u \alpha t } -1 )  + \tilde \lambda \int_0 ^ t e^x \int_x ^ {\infty} h (s) e^{-s} ds dx$. $\Psi$ is monotone: if for all $t \geq 0$, $h_1 (t) \geq  h_2 (t)$ then for all $t \geq 0$, $\Psi (h_1) (t) \geq \Psi (h_2) (t)$.   Since, for all $t \geq 0$,  $f_u(t)  \geq \Psi (f_u)(t) \geq 1$, we deduce by iteration that there exists a function $h$ such that $h = \Psi (h) \geq 1$. Solving $h = \Psi (h)$ is simple, taking twice the derivative, we get,
$
h'' - h' + \tilde \lambda h = - \epsilon e^{p \alpha t} .
$
Therefore, $h = a e^{\tilde \alpha t } + b e^{\tilde \beta t } - \epsilon ( u^2 \alpha^2 - u \alpha + \tilde \lambda ) ^{-1}e^{u \alpha t }  $ for some constant $a$ and $b$. From (\ref{eq:deftilde2}) the leading term as $t$ goes to infinity is equal to $- \epsilon ( u^2 \alpha^2 - u \alpha + \tilde \lambda ) ^{-1}e^{u \alpha t } $. However from (\ref{eq:defeps}), $- \epsilon ( u^2 \alpha^2 - u \alpha + \tilde \lambda ) ^{-1}  < 0$ and it contradicts the assumption that $h (t) \geq 1$ for all $t \geq 0$. Therefore we have proved that $u \leq \overline \gamma$.

\subsubsection{Proof of (\ref{eq:power2}).}
Let $f_u ^{(n)} (t) =\EE [\min (Y(t),n)^u]$, we have the following lemma.
\begin{lemma}
\label{le:powerupp}
There exists a constant $C > 0$ such that for all $t \geq 0$:
$$
f_u ^{(n)} (t)   \leq   C  e^{u \alpha t }   + \lambda \int_0 ^ t e^x \int_x ^ {\infty} f_u ^{(n)} (s) e^{-s} ds dx.
$$
\end{lemma}

The statement (\ref{eq:power2}) is a direct consequence of Lemmas  \ref{le:powerupp2} and \ref{le:powerupp}. Indeed, note that $f_u ^{(n)} \leq n^u$, thus by Lemma \ref{le:powerupp2}, for all $t \geq 0$, $f_u ^{(n)} (t)\leq C_1 e^{u \alpha t}$ for some positive constant $C_1$ independent of $n$. From the Monotone Convergence Theorem, we deduce that, for all $t \geq 0$, $f_u (t) \leq C_1 e^{u \alpha t}$. It remains to prove Lemma \ref{le:powerupp}.

\noindent{{\em Proof of Lemma \ref{le:powerupp}.}
The lemma is already proved if $u$ is an integer in (\ref{eq:boundfpn}). The general case is a slight extension of the same argument. We write $u = p - 1 + v$ with $v \in (0,1)$ and $p \in \N^*$. We use the inequality, for all $y_i \geq 0$, $1 \leq i \leq N$,
$$
\left( \sum_{i=1} ^N y_i  \right) ^ u \leq  \sum_{ i = 1} ^N \sum_{k=0} ^{p-1} { p-1 \choose k } y_i^{k+v}\left( \sum_{j \neq i} ^N y_i  \right) ^ {p-k-1} $$
(which follows from the inequality $(\sum y_i) ^{v} \leq \sum y_i^{v}$). Then from (\ref{eq:RDE2}) we get the stochastic domination
\begin{eqnarray*}
(Y (t) - 1) ^u  &  \leq_{st} &  \sum_{ \xi_i \leq t } Y_i ( \xi_i
+ D_i )^u \\
& & \hspace{30pt} + \sum_{ \xi_i \leq t }   \sum_{k=0} ^{p-2} {
p-1 \choose k } Y_i(\xi_i + D_i)^{k+v} \left( \sum_{\xi_j \neq
\xi_i \leq t } Y_j ( \xi_j + D_j )  \right) ^ {p-k-1}
\end{eqnarray*} From Lemma \ref{le:rec}, there exists $C$ such
that for all $1 \leq k \leq p-1$, $f_k (t) \leq C e^{k\alpha t}$
and and $\int_0 ^ t \EE[ Y(x+D) ^ k] dx \leq C e^{k \alpha t}$.
Note also, by Jensen inequality, that for all $1 \leq k \leq p-2$,
$f_{k+v} (t) \leq f_{p-1} (t) ^{(k+v)/(p-1)} \leq C e
^{(k+v)\alpha t}$. The same argument (with $p$ replaced by $u$)
which led to (\ref{eq:boundfpn}) in the proof of Lemma
\ref{le:rec2} leads to the result. \ep

\subsection{Some comments on the birth-and-assassination process}

\subsubsection{Computation of higher moments}
\label{subsub:high}
It is probably hard to derive an expression for all moments of $N$, even if in the proof of Lemma \ref{le:rec2}, we have built an expression of the cumulants of $Y(t)$ by recursion. However, exact formulas become quickly very complicated. The third moment, computed by hand, gives
\begin{equation*}
\label{eq:EY3}
f_3(t)  =  3 \frac{3 \lambda - \alpha}{4 \lambda  - 3 \alpha}  e^{3\alpha t}  - 6 \frac{\lambda( 2 \lambda - \alpha)}{(3 \lambda -2 \alpha)^2} e^{2 \alpha t} +  \left(1+ 6 \frac{\lambda( 2 \lambda - \alpha)}{(3 \lambda -2 \alpha)^2}  - 3  \frac{3 \lambda - \alpha}{4 \lambda  - 3 \alpha} \right) e^{\alpha t} .
\end{equation*}
Since $N \stackrel{d}{=} Y (D)$, we obtain,
\begin{eqnarray*}
\EE N ^3 & = & 6 \frac{(3 \lambda - \alpha) \alpha}{(4 \lambda  - 3 \alpha)(1 - \alpha - 3 \lambda)} - 6\frac{\lambda( 2 \lambda - \alpha)\alpha }{(3 \lambda -2 \alpha)^2(1 - \alpha - 2 \lambda)} + \frac{1}{1- \alpha}.
\end{eqnarray*}

\subsubsection{Integral equation of the Laplace transform}

It is also possible to derive an integral equation for the Laplace transform of $Y(t)$: $L_\theta (t) = \EE \exp(-\theta Y(t))$, with $\theta >0$. Indeed, using RDE (\ref{eq:RDE2}) and the exponential formula (\ref{eq:expformula}),
\begin{eqnarray*}
L_\theta (t) & = &e ^{-\theta} \exp \left( \lambda \int_0 ^t (\EE L_\theta (x+D) -1 ) dx\right) \\
& = & e ^{-\theta} \exp \left( \lambda \int_0 ^t  e^x \int_x ^\infty (L_\theta (s) -1 )e^{-s} ds dx\right).
\end{eqnarray*}
Taking twice the derivative, we deduce that, for all $\theta >0$, $L_\theta$ solves the differential equation:
$$
x'' x - {x' }^2 - x'x + \lambda x^2 (x -1) = 0.
$$
We have not been able to use fruitfully this non-linear differential equation.

\subsubsection{Probability of extinction}

If $\lambda > 1/4$ from Corollary \ref{cor:stab}, the probability of extinction of $\B$ is strictly less than $1$. It would be very interesting to have an asymptotic formula for this probability as $\lambda$ get close to $1/4$ and compare it with the Galton-Watson process. To this end, we define $\pi(t)$ as the probability of extinction of $\B$ given than the root cannot die before $t$. With the notation of Equation (\ref{eq:RDE2}), $\pi(t)$ satisfies
$$
\pi(t) = \EE \prod_{ i : \xi_i \leq t + D} \pi(t+D-\xi_i) = \EE \prod_{ i : \xi_i \leq t + D} \pi(\xi_i),
$$
Using the exponential formula (\ref{eq:expformula}), we find that the function $\pi$ solves the integral equation:
$$
\pi(t) = e^t \int_{t}^{\infty} \exp\left(-(\lambda+1) s + \lambda \int_o ^ s \pi(x) dx \right)ds.
$$
After a quick calculation, we deduce that $\pi$ is solution of  the second order non-linear differential equation
$$
\frac{ x' - x''} { x - x' } = \lambda ( x -1).
$$
Unfortunately, we have not been able to get any result on the
function $\pi(t)$ from this differential equation.

\section{Rumor scotching in a complete network}
\label{sec:rumor}
\subsection{Definition and result} \label{subsec:rumor}

\begin{figure}[htb]
\begin{center}
\includegraphics[angle=0,height = 4cm]{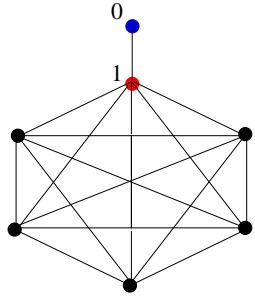}
\caption{The graph $G_6$.} \label{fig:k6}
\end{center}\end{figure}

We consider the rumor scotching process on the graph $G_n$ on
$\{0,\cdots,n\}$ obtained by adding on the complete graph on
$\{1,\cdots,n\}$ the edge $(0,1)$, see Figure \ref{fig:k6}. Let
$\cP_n$ be the set of subsets of $\{0,\cdots,n\}$. With the
notation in introduction, the rumor scotching process on $G_n$ is
the Markov process on $\cX_n = (\cP_n \times \{S,I,R\})^n$ with
generator, for $X = (A_i,s_i)_{0 \leq i \leq n}$,
\begin{align*}
& K(X, X + E_{ij})  =  \lambda n^{-1} \ind ( s_i = I) \ind ( s_j \neq R ), \\
& K(X, X - E_{j})  =  \ind( s_j = I  ) \left(\sum_{i  =1} ^ n \ind ( i \in A_j)  \right) ,
\end{align*}
and all other transitions have rate $0$. At time $0$, the initial
state is $X(0) = (X_i(0))_{0 \leq i \leq n}$ with $X_0 (0) = (\emptyset,R)$,  $X_1 (0) =
(\{0\},I)$ and for $i \geq 2$, $X_i (0) = (\emptyset,S)$.

With this initial condition, the  process describes the propagation of a rumor started from vertex $1$ at time $0$. After an exponential time, vertex $1$ learns that the rumor is false and starts to scotch the rumor to the vertices it had previously informed. This process is a Markov process on a finite set with as absorbing states, all states without $I$-vertices. We define $N_n$ as the total number of recovered vertices when the process stops evolving. We also define $Y_n (t)$ as the distribution $N_n$ given that vertex $1$ is recovered at time $t$. We have the following

\begin{theorem}\label{th:rq}
\begin{itemize}
\item[(i)]
If $0 < \lambda \leq 1/4$ and $t \geq 0$, as $n$ goes to infinity, $N_n$ and $Y_n (t)$ converge weakly respectively to $N$ and $Y (t)$ in the birth-and-assassination process of intensity $\lambda$.
\item[(ii)]
If $\lambda > 1/4$, there exists $\delta > 0$ such that $$\liminf_n \PP_{\lambda} ( N_n \geq \delta  n) > 0.$$
\end{itemize}
\end{theorem}

The proof of Theorem \ref{th:rq} relies on the convergence of the rumor scotching process to the birth-and-assassination process, exactly as the classical SIR dynamics converges to a branching process as the size of the population goes to infinity.

\subsection{Proof of Theorem \ref{th:rq}}

\subsubsection{Proof of Theorem \ref{th:rq}(i)} \label{subsubsec:rqi}

The proof of Theorem \ref{th:rq} relies on an explicit contruction of the rumor scotching process. Let $(\xi^ {(n)}_{ij}), 1 \leq i < j \leq n,$ be a collection of independent exponential variables with parameter $\lambda n^{-1}$ and, for all $1 \leq i \leq  j$,  let $D_{i j}$ be an independent exponential variable with parameter $1$. We set $D_{ji} = D_{ij}$ and $\xi^ {(n)}_{ji} = \xi^{(n)}_{ij}$. A network being a graph with marks attached on edges, we define $\cK_n$ as the network on the complete graph of $\{1,\cdots,n\}$ where the mark attached on the edge $(ij)$ is the pair $( \xi_{ij}^{(n)}, D_{ij})$. Now, the rumor scotching process is built on the network $\cK_n$ by setting $\xi^{(n)}_{ij}$ as the time for the infected particle $i$ to infect the particle $j$ and $D_{ij}$ as the time for the recovered particle $i$ to recover the particle $j$ that it had previously infected.

The network $\cK_n$ has a local weak limit as $n$ goes to infinity
(see Aldous and Steele \cite{aldoussteele} for a definition of the
local weak convergence). This limit network of $\cK_n$ is $\cK$,
the Poisson weighted infinite tree (PWIT) which is described as
follows. The root vertex, say $\emptyset$, has an infinite number
of children indexed by integers. The marks associated to the edges
from the root to the children are $(\xi_i, D_i)_{i \geq 1}$ where
$\{\xi_i\}_{i \geq 1}$ is the realization of a Poisson  process of
intensity $\lambda$ on $\R_+$ and $(D_i)_{i \geq 1}$  is a
sequence of independent exponential variables with parameter $1$.
Now recursively, for each vertex $i \geq 1$ we associate an
infinite number of children denoted by $(i,1), (i,2), \cdots$ and
the marks on the edges from $i$ to its children are obtained from
the realization of an independent Poisson  process of intensity
$\lambda$ on $\R_+$ and a sequence of independent exponential
variables with parameter $1$. This procedure is continued for all
generations. Theorem 4.1 in \cite{aldoussteele} implies the local
weak convergence of $\cK_n$ to $\cK$ (for a proof see Section 3 in
Aldous \cite{aldous92}).

Now notice that the birth-and-assissination process is the rumor scotching process on $\cK$ with initial condition: all vertices susceptible apart from the root which is infected and will be restored after an exponential time with mean $1$.

For $s > 0$ and $\ell \in \N$, let $\cK_n [s,\ell]$ be the network spanned by the set of vertices $j \in \{1,\cdots,n\}$ such that there exists a sequence $(i_1,\cdots,i_k)$ with $i_1 = 1$, $i_k = j$, $k \leq \ell$ and $\max (\xi^ {(n)}_{i_{1} i_{2}}, \cdots, \xi^ {(n)}_{i_{k-1} i_{k}}) \leq s$. If $\tau_n$ is the time elapsed before an absorbing state is reached, we get that $\ind(\tau_n \leq s) \ind(N_n \leq \ell)$ is measurable with respect to $\cK_n [s,\ell]$. From Theorem 4.1 in \cite{aldoussteele}, we deduce that $\ind(\tau_n \leq s) \ind(N_n \leq \ell)$ converges in distribution to $\ind(\tau \leq s) \ind(N \leq \ell)$ where $\tau$ is the time elapsed before all particles die in the birth-and-assassination process. If $0 < \lambda < 1/4$, $\tau$ is almost surely finite and we deduce the statement (i).

\subsubsection{Proof of Theorem \ref{th:rq}(ii)}

In order to prove part (ii) we couple the birth-and-assassination process and the rumor scotching process. We use the above notation and  build the rumor scotching process on the network $\cK_n$. If $X = ( (A_i,s_i)_{0 \leq i \leq n} )  \in \cX_n$, we define $I(X) = \{ 1 \leq i \leq n : s_i = I\}$ and $S(X) = \{ 1 \leq i \leq n : s_i = S\}$.

Let $X = X_n (u) \in \cX_n$ be the state of the rumor scotching process at time $u \geq 0$.  Let $i \in I(X)$, we reorder the variables $(\xi_{ij}^{(n)})_{ j \in S(X)}$ in non-decreasing  order: $\xi^{(n)}_{i j_1}\leq \cdots \leq \xi^{(n)}_{i j_{|S(X)|}}$. Define $\xi^{(n)}_{i j_{0}} = 0$, from the memoryless property of the exponential variable, for $1\leq k \leq |S(X)|$,  $\xi^{(n)}_{i j_{k}} - \xi^{(n)}_{i j_{k-1}}$  is an exponential variable with parameter $\lambda (|S(X)|-k+1)/n$  independent of $(\xi^{(n)}_{i j_{\ell}} - \xi^{(n)}_{i j_{\ell-1}}$, $\ell <  k)$. Therefore, for all $1 \leq k \leq |S(X)|$, the vector $(\xi^{(n)}_{i j_1}, \cdots , \xi^{(n)}_{i j_k})$ is stochastically dominated component-wise  by the vector $(\xi_{1}, \cdots , \xi_{k})$ where $ \{ \xi_j \}_{j \geq 1}$ is a Poisson  process of intensity $\lambda  ( |S(X)| - k +1)/n$ on $\R_+$ (i.e. for all $0 \leq t_1 \leq \cdots \leq t_{k}$, $ \PP (\xi^{(n)}_{i j_1} \geq t_1 , \cdots , \xi^{(n)}_{i j_{k}} \geq t_k) \leq \PP (\xi_{1} \geq t_1 , \cdots , \xi_{k} \geq t_{k})$).  In particular if $|S(X)| \geq (1 - \delta ) n$, with $0 < \delta < 1/2$, then $(\xi^{(n)}_{i 1}, \cdots , \xi^{(n)}_{i \lfloor n \delta \rfloor})$ is stochastically dominated component-wise by the first $\lfloor n \delta \rfloor$ arrival times of a Poisson process of intensity $\lambda (1 -2\delta)$.

Now, let $\delta > 0$ such that $\lambda' = \lambda(1 - 2\delta) > 1/4$. We define $S^{(n)}_u,I^{(n)}_u,R^{(n)}_u,$ as the number of $S,I,R$-particles at time $u \geq 0$ in $\cK_n$, and $I'_u$ as the number of particles "at risk" at time $u$ in the birth-and-assassination process with intensity $\lambda'$. Let $\tau_n = \inf \{ u \geq 0 : S^{(n)}_u \leq (1-\delta) n \}$. Note that if $0 \leq u \leq \tau_n$ then any $I$-particle has infected less than $\lfloor \delta n \rfloor$ $S$-particles. From what precedes, we get
$$
S^{(n)}_u \ind (u  \leq \tau_n)  \leq_{st}  n - I'_u.
$$
So that $S^{(n)}_u \leq_{st} \max( n - I'_u, (1-\delta) n) $. In particular, since $N_n \geq \sup_{u \geq 0 }  ( n - S^{(n)}_u ) $,  we get
$$
\PP_\lambda ( N_n \geq \delta n ) \geq \PP_{\lambda'} ( \limsup_{u \to \infty} I'_u = \infty).
$$
Finally, it is proved in \cite{aldouskrebs} that if $\lambda' > 1/4$ then  $\PP_{\lambda'} ( \limsup_{u \to \infty} I'_u = \infty) >0$.

\section*{Acknowledgement}
I am grateful to David Aldous for introducing me to the birth-and-assassination process and for his support on this work. I am indebted to David Windisch for pointing a mistake in the previous proof of Theorem \ref{th:N}(ii)  at $\lambda =1/4$.

\bibliographystyle{abbrv}
\bibliography{../../bib}

\end{document}